\documentclass[11pt]{amsart}
\usepackage[a4paper,margin=28mm]{geometry}
\usepackage[T1]{fontenc}
\usepackage{lmodern,amsmath,amssymb,mathtools,microtype}
\usepackage[hidelinks]{hyperref}
\hypersetup{pdftitle={Inverse problem in difference Galois theory},pdfsubject={Algebraic realization and rational difference-invariant lines}}
\newtheorem{theorem}{Theorem}[section]
\newtheorem{lemma}[theorem]{Lemma}

\theoremstyle{definition}
\newtheorem{definition}[theorem]{Definition}
\newtheorem*{example}{Example}
\theoremstyle{remark}
\newtheorem{remark}[theorem]{Remark}
\newcommand{\GL}{\mathrm{GL}}
\newcommand{\SL}{\mathrm{SL}}
\newcommand{\Ga}{\mathbb G_a}
\newcommand{\Gm}{\mathbb G_m}
\newcommand{\Q}{\mathbb Q}
\newcommand{\Z}{\mathbb Z}
\newcommand{\Aut}{\operatorname{Aut}}
\newcommand{\Stab}{\operatorname{Stab}}
\newcommand{\Span}{\operatorname{span}}
\newcommand{\Lie}{\operatorname{Lie}}
\newcommand{\ord}{\operatorname{ord}}

\newcommand{\im}{\operatorname{im}}
\newcommand{\diag}{\operatorname{diag}}
\newcommand{\PV}{Picard--Vessiot}
\numberwithin{equation}{section}
\title{Inverse problem in difference Galois theory}
\author{Ruyong Feng}
\address{State Key Laboratory of Mathematical Sciences,
Academy of Mathematics and Systems Science, Chinese Academy of Sciences\\
School of Mathematical Sciences, University of Chinese Academy of Sciences}
\email{ryfeng@amss.ac.cn}
\date{}

\begin{document}
\begin{abstract}
Let $C$ be an algebraically closed field of characteristic zero, and let $K=C(x)$ with $\sigma(x)=x+1$. 
We show that every closed subgroup $G\subseteq\GL_n(C)$ with cyclic component group $G/G^\circ$ can be realized as the \PV{} group 
of a system $\sigma Y=A(x)Y$. Here $G^\circ$ denotes the identity component of $G$. 
Moreover, we can choose the coefficient matrix to have the formal expansion $A(x)=h+O(x^{-2})$ at infinity, 
where $h\in G(C)$ is semisimple. This gives an affirmative answer to the conjecture of van der Put and Singer in their monograph
\emph{Galois Theory of Difference Equations} (Lecture Notes in Mathematics, vol.~1666). 
The original results in this paper were entirely generated by ChatGPT 6 and the MechMath Agent Team (MMAT) developed 
by the Key Laboratory of Mathematics Mechanization, AMSS, CAS. The author reviewed and simplified 
the initial proofs and verified the correctness of the results.
\end{abstract}
\maketitle

\section{Introduction}\label{sec_original_1}
The inverse problem is one of the fundamental problems in Galois theory. It asks which groups can occur as Galois groups over a fixed base field. For linear differential and difference equations, \PV{} theory associates a linear algebraic group with each system. The corresponding inverse problem is to realize prescribed linear algebraic groups by such systems. The differential inverse problem over the rational function field $C(x)$ has been completely solved. For the shift operator $\sigma(x)=x+1$, however, the corresponding inverse problem had remained open in general, although positive results were known for several classes of groups.

\subsection{The differential inverse problem}\label{sub_original_1_1}
In his address to the 1966 International Congress of Mathematicians in Moscow, Kolchin \cite{Kolchin} called attention to the differential inverse problem: given a linear algebraic group, find a linear differential equation over a prescribed differential field having that group as its Galois group.

Early algebraic results of Bia\l ynicki-Birula and Kovacic treated connected nilpotent and connected solvable groups, respectively. Kovacic also developed reductions toward the reductive case. These developments and their relation to later constructions are reviewed in \cite[Introduction]{MS96}. Over $\mathbb C(x)$, with derivation $d/dx$, Tretkoff and Tretkoff \cite{TT} proved in 1979 that every complex linear algebraic group is a differential Galois group. Their analytic method uses the Riemann--Hilbert problem to realize a finitely generated Zariski-dense subgroup as monodromy. The density theorem for regular singular equations then identifies the resulting differential Galois group.

For an arbitrary algebraically closed constant field $C$ of characteristic zero, Mitschi and Singer \cite{MS96} gave a constructive, purely algebraic realization of every connected linear algebraic group. Their theorem applies more generally to differential fields of finite, nonzero transcendence degree over $C$ with constant field $C$. They subsequently treated solvable-by-finite groups over $C(x)$ \cite{MS02}. Hartmann \cite{Hartmann} completed the inverse problem over $C(x)$ for arbitrary linear algebraic groups, using algebraic methods.

More recent work has strengthened this existence theorem in several directions. Feng and Wibmer \cite{FW} established specialization results for differential Galois groups and proved Matzat's conjecture. This conjecture states that the absolute differential Galois group of a one-variable function field $F$ over $C$, equipped with a nonzero $C$-derivation, is a free proalgebraic group on a set of cardinality $|F|$. Together with \cite[Corollary 3.8(b)]{BHHW}, this result implies that, for every nontrivial linear algebraic group $G$ over $C$, there are exactly $|C|$ differential $C(x)$-isomorphism classes of \PV{} extensions $L/C(x)$ with differential Galois group isomorphic to $G$. Serafini and Wibmer \cite{SW} also proved that every linear algebraic group over $C$ can be realized as the differential Galois group of a regular singular differential equation over $C(x)$.

\subsection{The difference inverse problem}\label{sub_original_1_2}
For the difference field $(C(x),\sigma)$, where $C$ is algebraically closed of characteristic zero and $\sigma(x)=x+1$, van der Put and Singer \cite[Proposition 1.20]{PS} proved that $G/G^\circ$ is necessarily cyclic whenever $G$ is a \PV{} group. They stated the converse as a conjecture in Chapter~3, page~35 of \cite{PS}: every linear algebraic group $G$ for which $G/G^\circ$ is cyclic should occur as a difference Galois group over $C(x)$.

In their monograph \cite{PS}, van der Put and Singer laid the foundations for solving the difference inverse problem.
In particular, they introduced the coefficient group $T_A$ associated with a system $\sigma Y=AY$ 
(see Definition~\ref{result_original_2_1}). They showed that the difference Galois group of this system is 
conjugate to a subgroup of $T_A$ \cite[Proposition 3.2(2)]{PS}\footnote{The second assertion of Proposition~3.2(2) 
is false in general. Counterexamples can be found at \url{http://www.mmrc.iss.ac.cn/~ryfeng/papers/counterexamples.pdf} or at \url{https://singer.math.ncsu.edu/papers/Inverse_error.pdf}.}. Furthermore, they gave a criterion, in terms of rational changes of variables, for the coefficient group $T_A$ to be the difference Galois group \cite[Proposition 1.21]{PS}. Their algebraic constructions realize $\Ga^m$ and tori over $C$ as difference Galois groups over $C(x)$ \cite[Lemma 3.4(2)--(3)]{PS}. Over $\mathbb C(x)$, the analytic theory developed in Chapter~8 of \cite{PS} also realizes every connected linear algebraic group and every cyclic extension of a torus as a difference Galois group. The respective results are Corollary~8.6(3) and Lemma~8.12 of \cite{PS}.

Feng \cite{Feng} subsequently developed specialization results for difference Galois groups. These results transfer the above analytic realizations over $\mathbb C(x)$ to the corresponding realizations over $C(x)$ for an arbitrary algebraically closed field $C$ of characteristic zero \cite[Theorem 6.3 and Corollary 6.4]{Feng}. Thus every connected linear algebraic group over $C$, and every cyclic extension of a torus over $C$, can be realized as a difference Galois group over $C(x)$.

There is an analogous inverse Galois problem for the $q$-shift $\sigma_q(f(x))=f(qx)$. We do not discuss this problem in detail here and refer interested readers to Chapter~12 of the monograph of van der Put and Singer \cite{PS} and the article by Ramis and Sauloy \cite{RS}.

\subsection{Main result and conventions}\label{sub_original_1_3}
Throughout this paper, $C$ is algebraically closed of characteristic zero, $K=C(x)$, and $\sigma(f(x))=f(x+1)$. We identify the prime field with $\Q\subset C$. The field of constants of $K$, denoted by $K^\sigma$, is $C$. By an algebraic group we mean a reduced closed subgroup of a general linear group over $C$. We write $G^\circ$ for the connected component of $G$ containing the identity. Such groups and their closed subgroups are smooth in characteristic zero. A rational representation of $G$ is a regular homomorphism of algebraic groups $\rho:G\to\GL(V)$ on a finite-dimensional $C$-vector space; it is not a partially defined rational map. We use the same symbol for the induced map on points over any commutative $C$-algebra. The main result of this paper is the following theorem.

\begin{theorem}\label{result_original_1_1}\label{main}
Let $G\subseteq\GL_n(C)$, $n\geq1$, be a closed algebraic subgroup. Assume that $G/G^\circ$ is cyclic. There exist $A(x)\in G(K)$, a \PV{} ring $R$ for $\sigma Y=AY$, and a fundamental matrix $Z\in G(R)$ such that
\[
\Aut_\sigma(R/K)\longrightarrow\GL_n(C),\qquad\tau\longmapsto Z^{-1}\tau(Z)
\]
has image exactly $G(C)$.
\end{theorem}
The theorem establishes the sufficient direction of van der Put and Singer's conjecture. As noted in Section~\ref{sub_original_1_2}, their monograph already provides many of the tools needed to solve the inverse problem. In particular, they proved that whenever $G/G^\circ$ is cyclic, there exists $A\in G(K)$ with $T_A=G$ \cite[Proposition 3.3(2)]{PS}. Our main contribution is to construct such an $A$ satisfying $A(x)=h+O(x^{-2})$, with $h$ semisimple, and the singularity condition in \eqref{hyp}, and to establish a rigidity result for rational difference-invariant lines under these conditions. This rigidity, together with $T_A=G$, rules out any rational change of variables $U=BY$, with $B\in G(K)$, for which the transformed coefficient matrix $\sigma(B)AB^{-1}$ belongs to $H(K)$ for some proper closed algebraic subgroup $H\subsetneq G$.

\subsection{Contributions}\label{contributions}
The original proofs\footnote{The original proofs generated by AI tools can be found at \url{http://www.mmrc.iss.ac.cn/~ryfeng/papers/originalproofs-AI.pdf}} in this paper were generated entirely by ChatGPT 6 and the MechMath Agent Team (MMAT), 
after the author supplied the necessary reference materials. 
The original proofs contain two main new ingredients. 
The first is the construction of a coefficient matrix $A(x)$ satisfying 
$A(x)=h+O(x^{-2})$ at infinity, where $h$ is a semisimple constant matrix, together with the required singularity 
condition. The second is a rigidity theorem for rational difference-invariant lines (see Theorem~\ref{rigidity}).

Using these ingredients, the original AI-generated proof used descent over a Picard-Vessiot ring.
 Suppose that the difference Galois group $H$ of $\sigma Y=AY$ is realized as a proper closed subgroup of $T_A=G$. 
 By Chevalley's theorem, there exist an algebraic representation $\rho: G\to\mathrm{GL}(V)$ and 
 a line $\ell\subset V$ such that
$
H=\operatorname{Stab}_G(\ell).
$
Let $R$ be a Picard-Vessiot ring for the system, let $Q(R)$ be its total ring of fractions, 
and choose a fundamental matrix $Z\in G(R)$ realizing this identification of the Galois group. 
Descent produces a rational difference-invariant line $W\subset V_K=K\otimes_C V$ satisfying
$$
Q(R)W=Q(R)\rho(Z)\ell.
$$
The rigidity theorem implies that $W=K\ell_0$ for a constant line $\ell_0\subset V$. Since
$A\in\operatorname{Stab}_G(\ell_0)(K)$ and $T_A=G$,
the line $\ell_0$ is $G$-stable. Consequently, $\rho(Z)^{-1}$ preserves $Q(R)\ell_0$, and the preceding equality gives
$
Q(R)\ell=Q(R)\ell_0.
$
Thus $\ell=\ell_0$, contradicting
$
\operatorname{Stab}_G(\ell)=H\subsetneq G.
$
This argument requires extending scalars to $Q(R)$, invoking properties of the Picard-Vessiot ring, and carrying out a separate descent step.

The author observed that the proof could be simplified by applying Proposition 1.21 of van der Put and Singer's monograph directly, 
with all line constructions carried out over $K$ (see the proof of Theorem~\ref{criterion}). Specifically, suppose that
$$
\sigma(B)AB^{-1}\in H(K)
$$
for a proper closed subgroup $H\subsetneq G$ and some $B\in G(K)$. Choose a Chevalley representation $\rho$ and a line $\ell=Cw$ 
with stabilizer $H$. Then
$
W=K\rho(B^{-1})w
$
is rational difference-invariant. Rigidity gives $W=K\ell_0$ for a constant line $\ell_0\subset V$, and the identity
$$
\Stab_G(\ell_0)_K
=B^{-1}H_KB
$$
shows that its stabilizer is proper. On the other hand,
$
A\in\Stab_G(\ell_0)(K)
$
forces this stabilizer to contain $T_A=G$, a contradiction. The Galois group is therefore $G$ by \cite[Proposition 1.21(3)]{PS}.
 This simplification eliminates both the separate descent step and the subsequent use of $\rho(Z)^{-1}$. 
 Under the guidance of the author, the proof of Lemma 3.2 was also simplified by making greater use
  of existing results in Milne's monograph.

\subsection{Organization}\label{sub_original_1_4}
In Section~\ref{sec_original_2} we recall the coefficient group and the finite singularity set, study their behavior under algebraic representations, and prove the rigidity theorem for rational difference-invariant lines. We then combine rigidity with Chevalley's theorem to exclude rational reduction to any proper subgroup and apply van der Put and Singer's criterion to identify the Galois group. In Section~\ref{sec_original_6} we construct a rational curve through generating points of the identity component, impose the required expansion at infinity, and use a semisimple lift and an arithmetic dilation to complete the proof of Theorem~\ref{main}. In Section~\ref{sec_original_8} we give explicit realizations of $\SL_2$, $\Gm$, a nonsplit extension of a torus and finite cyclic groups, and use counterexamples to explain the roles of the hypotheses. In Section~\ref{sec_original_10} we give an alternative construction of the generating curve, using a rational chart and Hermite interpolation to prescribe tangent vectors that span the Lie algebra of the identity component.

\section{Rigidity of rational difference-invariant lines}\label{sec_original_2}
A rational difference-invariant line need not have a constant generator. We prove that every such line in every algebraic representation of $G$ is the scalar extension of a $C$-line when $A(x)=h+O(x^{-2})$ with $h$ semisimple and $1\notin\Span_\Q\mathcal{P}(A)$, where $\mathcal{P}(A)$ is the set of finite poles of $A$ and $A^{-1}$; the proof uses valuations at finite points and a comparison of coefficients at infinity. We then use this rigidity to exclude rational reduction to a proper subgroup when $T_A=G$, and deduce that the difference Galois group is $G$.

\subsection{Finite singularities and change of representation}
\begin{definition}\label{result_original_2_1}
For $A\in\GL_n(K)$, let
\[
 \mathcal{P}(A)=\{s\in C:\text{some entry of $A$ or $A^{-1}$ has a pole at $s$}\}.
\]
Let $T_A$ be the smallest closed $C$-subgroup of $\GL_n$ for which $A\in T_A(K)$. Equivalently, $T_A$ is the algebraic subgroup generated by the values $A(a)$ for $a\in C\setminus \mathcal{P}(A)$.
\end{definition}
The definition and the description of $T_A$ by its generating values are given in \cite[Proposition 3.2(1)]{PS}. Note that the poles of $A^{-1}$ are indispensable in the definition of $\mathcal{P}(A)$.

\begin{lemma}[Functorial control of singularities and the first jet]\label{result_original_2_2}
Let $A\in G(K)$ and let $\rho:G\to\GL(V)$ be an algebraic representation. Then
\[ \mathcal{P}(\rho(A))\subseteq \mathcal{P}(A). \]
If $A(x)=h+O(x^{-2})$ with $h\in G(C)$, then
\[ \rho(A(x))=\rho(h)+O(x^{-2}). \]
If $h$ is semisimple, so is $\rho(h)$.
\end{lemma}
\begin{proof}
For $s\notin \mathcal{P}(A)$, put $\mathcal O_s=C[x]_{(x-s)}$. Both $A$ and $A^{-1}$ have entries in $\mathcal O_s$, 
so $A\in\GL_n(\mathcal O_s)$. The defining equations of $G$ vanish on $A$ in $K$ and hence in $\mathcal O_s$. 
It follows that $A\in G(\mathcal O_s)$. Applying the morphism $\rho$ gives $\rho(A)\in\GL(V)(\mathcal O_s)$.
Thus both $\rho(A)$ and its inverse are regular at $s$.

At infinity, set $u=x^{-1}$. Then $A$ and $h$ are $G(C[[u]])$-points with the same reductions in $G(C[[u]]/(u^2))$. Applying $\rho$ preserves this congruence, proving the second assertion. 
The third assertion follows from the functoriality of Jordan decomposition.
\end{proof}

\begin{lemma}[Local primitivity]\label{result_original_2_3}
\leavevmode
\begin{enumerate}
\item[\textup{(i)}] Let $\mathcal O$ be a discrete valuation ring with fraction field $F$ and valuation $\nu$. For $0\ne w=(w_i)\in F^m$, put $\nu(w)=\min_i\nu(w_i)$, with $\nu(0)=+\infty$. Then, for $a\in F^*$ and $B\in\GL_m(\mathcal O)$,
\[
 \nu(aw)=\nu(a)+\nu(w),\qquad\nu(Bw)=\nu(w).
\]
\item[\textup{(ii)}] Every $K$-line $L\subset K^m$ has a primitive polynomial generator $v\in C[x]^m$, meaning $L=Kv$ and $\gcd(v_1,\ldots,v_m)=1$. It is unique up to multiplication by an element of $C^*$. Writing $\nu_s=\ord_{x=s}$ and extending it to vectors as in \textup{(i)}, we have
\[
 \nu_s(v)=\nu_s(\sigma v)=0\qquad\text{for every }s\in C.
\]
\end{enumerate}
\end{lemma}
\begin{proof}
\textup{(i)} The scalar identity follows directly from the definition. Since the entries of $B$ lie in $\mathcal O$,
\[
 \nu(Bw)
 =\min_i\nu\left(\sum_j B_{ij}w_j\right)
 \geq\min_{i,j}\bigl(\nu(B_{ij})+\nu(w_j)\bigr)
 \geq\nu(w).
\]
Applying this inequality to $B^{-1}\in\GL_m(\mathcal O)$ gives
\[
 \nu(w)=\nu\bigl(B^{-1}(Bw)\bigr)\geq\nu(Bw)\geq\nu(w).
\]

\textup{(ii)} Choose $0\ne u\in L$ and $b\in C[x]\setminus\{0\}$ such that $bu\in C[x]^m$. If $g$ is the monic greatest common divisor of the coordinates of $bu$, then $v=g^{-1}bu$ is a primitive polynomial generator of $L$. Since $C[x]$ is a principal ideal domain, there exist $a_1,\ldots,a_m\in C[x]$ such that
\[
 \sum_{i=1}^m a_i v_i=1,
 \qquad
 \sum_{i=1}^m \sigma(a_i)\sigma(v_i)=1.
\]
After localization at $(x-s)$, these identities show that the coordinates of each of $v$ and $\sigma v$ generate $\mathcal O_s=C[x]_{(x-s)}$. Thus $\nu_s(v)=\nu_s(\sigma v)=0$.

Finally, if $\widetilde v=fv$ is another primitive polynomial generator, with $f\in K^*$, then, as fractional ideals of $C[x]$,
\[
 C[x]=\sum_i C[x]\widetilde v_i
      =f\sum_i C[x]v_i=fC[x].
\]
Hence $f\in C[x]^*=C^*$, proving uniqueness.
\end{proof}

\subsection{Rational difference-invariant lines}\label{sec_original_3}
Throughout this section, let $V$ be a finite-dimensional $C$-vector space equipped with an algebraic representation 
$\rho:G\to\mathrm{GL}(V)$. Put $V_K=K\otimes_C V$, and identify $V$ with $1\otimes V$. 
Extend $\sigma$ to $V_K$ by
\[
\sigma(f\otimes v)=\sigma(f)\otimes v.
\]
This defines a $\sigma$-semilinear automorphism whose fixed vectors are precisely the elements of $V$, which we call constant vectors.

\begin{definition}[Rational difference-invariant line]\label{def_rational_difference_invariant_line}
Let $A\in G(K)$ and let $\rho:G\to\GL(V)$ be an algebraic representation. On $V_K$, consider the $\sigma$-semilinear automorphism
\[
 \Phi_\rho=\rho(A)^{-1}\sigma.
\]
A one-dimensional $K$-subspace $L\subset V_K$ is called a \emph{rational difference-invariant line} if $\Phi_\rho(L)=L$. Equivalently, for a nonzero generator $v$ of $L$, there exists $c\in K^*$ such that
\begin{equation}\label{line}
 \sigma(v)=c\rho(A)v.
\end{equation}
\end{definition}
Note that the condition in the above definition is independent of the chosen generator, because replacing $v$ by $fv$, where $f\in K^*$, replaces $c$ by $c\sigma(f)/f$.

\begin{theorem}[Rigidity of rational difference-invariant lines]\label{result_original_3_1}\label{rigidity}
Assume that $A\in G(K)$ satisfies
\begin{equation}\label{hyp}
 A(x)=h+O(x^{-2}),\qquad h\in G(C)\text{ semisimple},\qquad
 1\notin\Span_{\Q}\mathcal{P}(A).
\end{equation}
For every finite-dimensional algebraic representation $\rho$ of $G$, each rational difference-invariant $K$-line in $K\otimes_C V$ is the scalar extension of a $C$-line.
\end{theorem}
\begin{proof}
Fix a $C$-basis of $V$ and identify $V_K$ with $K^m$. Let $L\subset V_K$ be a rational difference-invariant line. By Lemma~\ref{result_original_2_3}\textup{(ii)}, it has a primitive polynomial generator $v\in C[x]^m$ satisfying
\begin{equation}\label{prim}
 \nu_s(v)=\nu_s(\sigma v)=0\qquad(s\in C).
\end{equation}
By Definition~\ref{def_rational_difference_invariant_line}, we have $\sigma(v)=c\rho(A)v$ for some $c\in K^*$. Put $B=\rho(A)$ and $\mathcal P=\mathcal P(A)$. For $s\notin\mathcal P$, Lemma~\ref{result_original_2_2} gives $B\in\GL_m(\mathcal O_s)$. Lemma~\ref{result_original_2_3}\textup{(i)} and \eqref{prim} therefore yield
\begin{equation}
 0=\nu_s(\sigma v)=\nu_s(cBv)
  =\nu_s(c)+\nu_s(v)=\nu_s(c).
\end{equation}
Thus $c$ has no finite zero or pole outside $\mathcal{P}$. Since $C$ is algebraically closed, we can write
\begin{equation}\label{factorc}
 c(x)=\gamma\prod_{s\in \mathcal{P}}(x-s)^{k_s},\qquad\gamma\in C^*,\quad k_s\in\Z.
\end{equation}

We now work with formal Laurent expansions at infinity, in $C((x^{-1}))$. Let $d=\max_i\deg v_i\geq0$ and write
\[
 v(x)=x^d\bigl(v_d+x^{-1}v_{d-1}+O(x^{-2})\bigr),\qquad v_d\ne0.
\]
If $d=0$, put $v_{-1}=0$. Set $e=\sum_s k_s$. Since $c(x)=\gamma x^e(1+O(x^{-1}))$ and $\rho(A(x))=\rho(h)+O(x^{-2})$, the highest powers of $x$ in the Laurent expansions of the two sides of \eqref{line} are respectively $d$ and $d+e$. The leading coefficient on the right is $\gamma\rho(h)v_d\ne0$. Therefore
\begin{equation}\label{sumzero}
 e=\sum_s k_s=0.
\end{equation}
Let $q=\sum_s k_s s$ and $M=\gamma\rho(h)$. Factoring powers of $x$ out of \eqref{factorc} and expanding formally, also for negative $k_s$, gives
\[ c(x)=\gamma\bigl(1-qx^{-1}+O(x^{-2})\bigr). \]
We consequently have
\begin{align*}
 v(x+1)&=x^d\bigl(v_d+x^{-1}(dv_d+v_{d-1})+O(x^{-2})\bigr),\\
 c\rho(A)v&=x^d\bigl(Mv_d+x^{-1}(Mv_{d-1}-qMv_d)+O(x^{-2})\bigr).
\end{align*}
Comparing coefficients yields
\begin{equation}\label{compare}
 Mv_d=v_d,\qquad(d+q)v_d=(M-I)v_{d-1}.
\end{equation}
By Lemma~\ref{result_original_2_2}, $M$ is semisimple. Its $1$-eigenspace has zero intersection with $\im(M-I)$. The two sides of the second equality in \eqref{compare} lie in these respective spaces. Since $v_d\ne0$, we get
\begin{equation}\label{degree}
 d+\sum_{s\in \mathcal{P}}k_s s=0.
\end{equation}
If $d>0$, division by the nonzero integer $d$ would put $1$ in $\Span_\Q \mathcal{P}$. Hence $d=0$, and $v$ is a nonzero constant vector.
\end{proof}

\begin{example}[A rational difference-invariant line with no constant generator]
The conclusion of Theorem~\ref{result_original_3_1} need not hold for a general coefficient matrix. Take $G=\GL_2$ with its natural representation on $V=C^2$, and put
\[
 A=\begin{pmatrix}1&1\\0&1\end{pmatrix},\qquad
 v(x)=\begin{pmatrix}x\\1\end{pmatrix},\qquad L=Kv\subset K\otimes_C V=K^2.
\]
Then
\[
 \sigma(v)=\begin{pmatrix}x+1\\1\end{pmatrix}=Av,
 \qquad\Phi_\rho(v)=A^{-1}\sigma(v)=v.
\]
Consequently, $\Phi_\rho(fv)=\sigma(f)v$ for every $f\in K$, and the surjectivity of $\sigma:K\to K$ gives $\Phi_\rho(L)=L$.
Nevertheless $L\cap V=\{0\}$. Indeed, if $fv=(a,b)^T$ with $f\in K$ and $a,b\in C$, comparison of coordinates gives $f=b$ and $bx=a$. Since $x$ is transcendental over $C$, these equalities imply $b=a=f=0$. If $L$ were the scalar extension $K\otimes_C L_0$ of a $C$-line $L_0\subseteq V$, any nonzero vector in $L_0$ would belong to $L\cap V$, a contradiction. Thus no rescaling of $v$ by an element of $K^*$ produces a constant vector.

Here $\mathcal{P}(A)=\varnothing$ and $A(x)=h+O(x^{-2})$ with $h=A$, but $h$ is not semisimple, so the hypotheses of Theorem~\ref{result_original_3_1} are not satisfied. The \PV{} group of this same system is discussed in Section~\ref{sub_original_9_2}.
\end{example}

\begin{remark}[Uniformity and the exact conclusion]\label{result_original_3_2}
The hypotheses are imposed on $A$ once, independently of $\rho$. In every representation the finite singularities lie in $\mathcal{P}(A)$ and the first-order term at infinity vanishes. What becomes constant is the line, not necessarily $c$; the integers $k_s$ need not vanish separately. For a constant generator $v_0$, \eqref{line} says that $\rho(A)$ preserves $Kv_0$, or $A\in\Stab_G(Cv_0)(K)$. It does not say that $\rho(A)$ fixes $v_0$.
\end{remark}
\begin{remark}[A translation-invariant weakening]\label{result_original_3_3}
If $\mathcal{P}$ is nonempty, the last condition in \eqref{hyp} may be replaced by
\[ 1\notin\Span_\Q\{s-s_0:s\in \mathcal{P}\},\qquad s_0\in \mathcal{P}. \]
Indeed, \eqref{sumzero} gives $\sum_s k_s s=\sum_s k_s(s-s_0)$. This condition is independent of $s_0$. Equivalently,
\[
 \left\{\sum_s k_s s:k_s\in\Z,\ \sum_s k_s=0\right\}\cap\Z=\{0\}.
\]
One implication follows by clearing denominators, and the other by dividing by a nonzero integer. The stronger condition \eqref{hyp} is sufficient for the construction below.
\end{remark}

\subsection{Rational reduction and the coefficient-group criterion}\label{sub_rational_reduction}
The change of variables $U=BY$, with $B\in G(K)$, transforms $\sigma Y=AY$ into
\[
 \sigma U=\sigma(B)AB^{-1}U.
\]
We shall show that, under the rigidity hypotheses, the transformed matrix cannot belong to a proper subgroup of $G=T_A$. The Galois group will then be determined by \cite[Proposition 1.21(3)]{PS}.

\begin{theorem}[Coefficient-group criterion]\label{result_original_5_2}\label{criterion}
Suppose that $A\in G(K)$ satisfies \eqref{hyp} and $T_A=G$. For every proper closed $C$-subgroup $H\subsetneq G$ and every $B\in G(K)$, we have
\[
 \sigma(B)AB^{-1}\notin H(K).
\]
Consequently, the difference Galois group of $\sigma Y=AY$ is $G$.
\end{theorem}
\begin{proof}
Suppose that $D=\sigma(B)AB^{-1}\in H(K)$ for some $B\in G(K)$ and some proper closed subgroup $H\subsetneq G$. By Chevalley's theorem \cite[Theorem 4.27]{Milne}, there are an algebraic representation $\rho:G\to\GL(V)$ and a line $\ell=Cw$ such that $H=\Stab_G(\ell)$. Write $\rho(D)w=aw$ with $a\in K^*$, and put $v=\rho(B^{-1})w$. Since $\sigma(B)^{-1}=AB^{-1}D^{-1}$, we have
\[
 \sigma(v)=\rho\bigl(\sigma(B)^{-1}\bigr)w
 =\rho(A)\rho(B^{-1})\rho(D^{-1})w
 =a^{-1}\rho(A)v.
\]
Thus $Kv$ is a rational difference-invariant line. Theorem~\ref{rigidity} gives $Kv=K\ell_0$ for a constant line $\ell_0\subset V$.

Let $H_1=\Stab_G(\ell_0)$. For every commutative $K$-algebra $S$, the equality $K\ell_0=\rho(B^{-1})(K\ell)$ gives $S\ell_0=\rho(B^{-1})(S\ell)$. Thus, for $g\in G(S)$,
\[
 \rho(g)(S\ell_0)=S\ell_0
 \quad\Longleftrightarrow\quad
 \rho(BgB^{-1})(S\ell)=S\ell.
\]
Since $H=\Stab_G(\ell)$, this says that $g\in H_1(S)$ if and only if $BgB^{-1}\in H(S)$. The equality of these subgroup functors proves
\[
 (H_1)_K=B^{-1}H_KB.
\]
Hence $H_1$ is a proper closed subgroup of $G$. Since $\ell_0$ is constant, $\sigma(K\ell_0)=K\ell_0$. Difference invariance therefore gives $\rho(A)(K\ell_0)=K\ell_0$, or $A\in H_1(K)$. It follows that
\[
 G=T_A\subseteq H_1\subsetneq G,
\]
a contradiction. The Galois-group assertion follows from \cite[Proposition 1.21(3)]{PS}.
\end{proof}

\section{Construction of realizing systems}\label{sec_original_6}
In this section we construct the coefficient matrices needed to prove Theorem~\ref{main}. We first interpolate through finitely many generating points of the identity component while making the first-order term at infinity vanish. We then choose a semisimple representative of a component-group generator and dilate the finite singularities so that the coefficient-group criterion applies.

\subsection{Unipotent elements and finitely many one-parameter subgroups}\label{sub_original_6_1}
Both lemmas in this subsection follow from classical structure theory. We give elementary proofs using polynomial exponential maps, Jordan decomposition, and the containment of semisimple elements in tori \cite[Theorem 9.18 and Corollary 17.36]{Milne}. The subsequent use of finitely many generating points follows the strategy of \cite[Proposition 3.3(2)]{PS}.

\begin{lemma}[A unipotent element lies in the identity component]\label{result_original_6_1}
If $u\in G(C)$ is unipotent, then $u\in G^\circ(C)$. More precisely, if $u\ne I$ and $N=\log(u)$, the polynomial homomorphism
\[ \Ga\longrightarrow\GL_n,\qquad t\longmapsto\exp(tN) \]
factors through $G$ and maps $1$ to $u$.
\end{lemma}
\begin{proof}
The logarithm series for a unipotent matrix and the exponential series for a nilpotent matrix both terminate. 
Their formal identities give $\exp((s+t)N)=\exp(sN)\exp(tN)$ and $\exp(N)=u$. 
For every integer $m$, $\exp(mN)=u^m\in G(C)$. For any defining equation $f\in I(G)$, 
the function $f(\exp(tN))$ is a polynomial in $t$: since $\det(\exp(tN))=1$, 
the pullback of every regular function on $\mathrm{GL}_n$ is a polynomial in $t$.
 It vanishes at all integers, an infinite subset of $C$, and is therefore identically zero. The image is consequently contained in $G$. Since it is connected and contains $I$, it lies in $G^\circ$.

If $u\ne I$, choose an entry $N_{ij}\ne0$. The polynomial
\[
 \frac1{N_{ij}}\left(\sum_{k=1}^{n-1}\frac{(-1)^{k+1}}{k}(X-I)^k\right)_{ij}
\]
pulls back to $t$, since $\log(\exp(tN))=tN$. Thus the induced map $C[G]\to C[t]$ is surjective, showing that the parametrization is a closed immersion of $\Ga$. The case $u=I$ is immediate.
\end{proof}

\begin{lemma}[Finite one-parameter generation]\label{result_original_6_2}
Every nontrivial connected linear algebraic group $H$ over $C$ is generated as an algebraic group by finitely many closed subgroups isomorphic to $\Ga$ or $\Gm$.
\end{lemma}
\begin{proof}
Let $\mathcal U$ be the family of closed subgroups of $H$ isomorphic to $\Ga$ or $\Gm$. For $\mathcal F\subseteq\mathcal U$, denote by $P_{\mathcal F}$ the algebraic subgroup generated by $\mathcal F$, whose existence follows from \cite[Proposition 2.51]{Milne}.

For $g\in H(C)$, write $g=g_sg_u$ for its multiplicative Jordan decomposition in $H$ \cite[Theorem 9.18]{Milne}. By Lemma~\ref{result_original_6_1}, every nonidentity $g_u$ belongs to a closed subgroup isomorphic to $\Ga$. By \cite[Corollary 17.36]{Milne}, $g_s$ belongs to a torus in $H$. Since $C$ is algebraically closed, this torus is a product of copies of $\Gm$ \cite[Definition 12.14 and Theorem 12.18(d)]{Milne}. Thus every $g\in H(C)$ belongs to $P_{\mathcal U}(C)$, and \cite[Corollary 1.18]{Milne} gives
\begin{equation}\label{fullgen}
 P_{\mathcal U}=H.
\end{equation}

For a nonempty finite subfamily $\mathcal F=\{U_1,\ldots,U_r\}$, the multiplication morphism
\[
 U_1\times\cdots\times U_r\longrightarrow H,\qquad
 (u_1,\ldots,u_r)\longmapsto u_1\cdots u_r
\]
has a smooth, connected source. Its image is contained in $P_{\mathcal F}$ and contains every $U_i$, 
as is seen by setting all other coordinates equal to the identity. 
Hence its image generates $P_{\mathcal F}$. By \cite[Proposition 2.53 and Corollary 1.35]{Milne}, $P_{\mathcal F}$ is connected and irreducible. This also holds for $P_\varnothing=\{1_H\}$.

Choose a finite $\mathcal F_0\subseteq\mathcal U$ for which $\dim P_{\mathcal F_0}$ is maximal; such a choice exists because these dimensions are integers between $0$ and $\dim H$. For every $U\in\mathcal U$,
\[
 P_{\mathcal F_0}\subseteq P_{\mathcal F_0\cup\{U\}},\qquad
 \dim P_{\mathcal F_0}=\dim P_{\mathcal F_0\cup\{U\}}.
\]
The larger subgroup is irreducible, so the two subgroups are equal \cite[Appendix A.47]{Milne}. Hence $P_{\mathcal F_0}$ contains every $U\in\mathcal U$, and \eqref{fullgen} implies $P_{\mathcal F_0}=H$. The family $\mathcal F_0$ is nonempty because $H$ is nontrivial.
\end{proof}

Fix such a finite list, with embeddings
\[ u_i:\Ga\hookrightarrow H\quad\text{or}\quad\lambda_i:\Gm\hookrightarrow H,\qquad1\leq i\leq r. \]
For an additive subgroup, put $g_i=u_i(1)$; for a multiplicative subgroup, put $g_i=\lambda_i(2)$. The powers of $g_i$ are respectively $u_i(m)$, $m\in\Z$, or $\lambda_i(2^m)$, $m\in\Z$. These parameter sets are infinite in characteristic zero. Every infinite subset of $\mathbb A^1$ or $\Gm$ is Zariski dense, since a nonzero polynomial in one variable has finitely many roots. Therefore
\begin{equation}\label{pointgen}
 H=\overline{\langle g_1,\ldots,g_r\rangle}^{\mathrm{Zar}}.
\end{equation}
Only these finitely many points will be prescribed on a rational curve.

\subsection{Interpolation with a vanishing first jet}\label{sub_original_6_2}
The formulas here give the required second-order condition in the coefficient variable $x$, using elementary Lagrange interpolation on one-parameter subgroups. Their singularities can then be controlled independently of the matrix entries defining those subgroups. For $1\leq i\leq r$, define
\begin{equation}\label{interpolants}
 M=r+1,\qquad L_i(z)=\prod_{\substack{1\leq j\leq r\\j\ne i}}\frac{z-j}{i-j},
 \qquad p_i(z)=\frac{i^M L_i(z)}{z^M}.
\end{equation}
The empty product is understood to be $1$. Since distinct nonzero integers remain distinct in $C$, these formulas are valid. They satisfy
\begin{equation}
 p_i(j)=\delta_{ij},\qquad p_i(z)=O(z^{-2})\quad(z=\infty).
\end{equation}
The second assertion follows from $\deg L_i=r-1=M-2$. Put
\begin{equation}\label{curve}
 \begin{aligned}
 F_i(z)&=\begin{cases}u_i(p_i(z)),&\text{for an additive subgroup},\\
 \lambda_i(1+p_i(z)),&\text{for a multiplicative subgroup},\end{cases}\\
 B(z)&=F_1(z)\cdots F_r(z).
 \end{aligned}
\end{equation}
The rational function $1+p_i$ is nonzero because its value at infinity is $1$, so the multiplicative expression defines a point over $C(z)$. All factors are $H(C(z))$-points.

At $z=j$, every additive parameter is $0$ except possibly the $j$th, and every multiplicative parameter is $1$ except possibly the $j$th, which is $2$. All evaluations are defined. Hence
\begin{equation}\label{curvevalues}
 B(j)=g_j,\qquad B(z)=I+O(z^{-2}),\qquad B(\infty)=I.
\end{equation}
The order estimate holds because each regular group homomorphism maps its identity to $I$ and 
preserves congruence modulo $z^{-2}$. Multiplying finitely many matrices congruent to $I$ modulo $z^{-2}$ preserves 
that congruence. By \eqref{pointgen}, the values of $B$ generate $H$.

\subsection{A finite set of algebraic numbers containing the singularities}\label{sub_original_6_3}
The entries of $u_i(t)$ and $u_i(t)^{-1}=u_i(-t)$ are polynomials in $t$. Thus an additive factor in \eqref{curve} has no finite pole except possibly at $z=0$. The entries of $\lambda_i(t)$ and its inverse are Laurent polynomials, since $C[\Gm]=C[t,t^{-1}]$. Consequently, a multiplicative factor and its inverse have finite poles only among
\[ 0\quad\text{and the zeros of}\quad Q_i(z)=z^M+i^M L_i(z). \]
These $Q_i$ are monic polynomials in $\Q[z]$. Thus
\begin{equation}\label{Szero}
 \mathcal{P}_0=\{0\}\cup\bigcup_{i\ \mathrm{multiplicative}}
 \{\alpha\in\overline\Q:Q_i(\alpha)=0\}
 \quad\text{satisfies}\quad \mathcal{P}(B)\subseteq \mathcal{P}_0.
\end{equation}
Here we choose an algebraic closure of $\Q$ inside $C$. Products and inverse products create no pole outside the union of the factor singularities. Possible cancellations only shrink $\mathcal{P}(B)$. The points of $\mathcal{P}_0$ are algebraic even when the entries defining $H$ and its one-parameter subgroups are arbitrary constants in $C$.

\subsection{Arithmetic dilation and realization of the full group}\label{sec_original_7}
The dilation below is an elementary avoidance argument in a finite-dimensional $\Q$-vector space. 
Its role in the proof is to arrange the hypothesis of Theorem~\ref{result_original_3_1}. 

\begin{lemma}[Dilation of a finite singularity set]\label{result_original_7_1}
For every finite $\mathcal{P}_0\subset C$, there is $t\in C^*$ such that $1\notin\Span_\Q(t\mathcal{P}_0)$. If $\mathcal{P}_0$ lies in a number field $F$ of degree $D$, any root $t\in C$ of $T^{D+1}-2$ is such a choice.
\end{lemma}
\begin{proof}
The space $W=\Span_\Q \mathcal{P}_0$ is finite dimensional. The field $C$ has infinite dimension over $\Q$, since it contains algebraic numbers of arbitrarily large degree. Choose $a\in C\setminus W$ and put $t=a^{-1}$. Note that $a\ne0$, since $0\in W$. The relation $1\in tW$ would imply $a\in W$.

For the explicit assertion, Eisenstein's criterion gives $[\Q(t):\Q]=D+1$. Hence $t\notin F$ and $t^{-1}\notin F$. Since $\Span_\Q(t\mathcal{P}_0)\subseteq tF$, this space cannot contain $1$.
\end{proof}

We shall also use the following semisimple lift of a generator of the component group.

\begin{lemma}[A semisimple lift of a quotient generator]\label{result_original_7_2}
If $G/G^\circ$ is cyclic, there is a semisimple $h\in G(C)$ whose image generates that quotient. Moreover, $G(C)=\langle G^\circ(C),h\rangle$.
\end{lemma}
This is the group-theoretic observation in \cite[proof of Theorem 8.10, p.~91]{PS}; it applies over every algebraically closed field of characteristic zero.
No assertion that $h^m=I$ is used. In particular, the exact sequence $1\to G^\circ\to G\to G/G^\circ\to1$ need not split, and no finite-order lifting result is required.

\begin{proof}[Proof of Theorem~\ref{result_original_1_1}]
Put $H=G^\circ$. Suppose first that $H\ne1$. Choose $B$ from \eqref{curve}, $h$ from Lemma~\ref{result_original_7_2}, and $t$ from Lemma~\ref{result_original_7_1} for the set \eqref{Szero}. Define
\begin{equation}\label{realize}
 A(x)=hB(x/t).
\end{equation}
Then $A\in G(K)$ and
\[
 A(x)=h+O(x^{-2}),\qquad \mathcal{P}(A)\subseteq t\mathcal{P}_0,\qquad1\notin\Span_\Q \mathcal{P}(A).
\]
Left multiplication by the constant $h$ and right multiplication by its inverse introduce no singularities.

It remains to prove that $T_A=G$. The group $T_A$ is contained in $G$. For every equation $f\in I(T_A)$, the rational function $f(A(x))$ is zero. The matrix $A$ is regular and invertible at infinity, so evaluation there gives $f(h)=0$. Thus $h\in T_A(C)$. Also, using \eqref{curvevalues}, we obtain
\[ A(ti)=hg_i,\qquad g_i=h^{-1}A(ti)\in T_A(C). \]
These evaluations are defined and invertible, as checked above. Since the $g_i$ generate $H$ and the image of $h$ generates $G/G^\circ$, we have $T_A=G$. Theorem~\ref{result_original_5_2} gives the required realization.

If $H=1$, then $G$ is a finite cyclic group. Take its generator $h$ and put $A=h$. In characteristic zero a finite-order matrix is semisimple, because its minimal polynomial divides a separable polynomial $T^m-1$. Now $\mathcal{P}(A)=\varnothing$, $A=h+O(x^{-2})$, and $T_A=G$. The same criterion applies, including to the trivial group.
\end{proof}

\section{Examples and counterexamples}\label{sec_original_8}\label{sec_original_9}
In this section we give explicit realizations of several algebraic groups and examine the hypotheses of the coefficient-group criterion. The counterexamples show why the finite singularities, semisimplicity and the first-order term at infinity must be controlled.

\subsection{A realization of $\mathrm{SL}_2$ with a single finite singularity}\label{sub_original_8_1}
Let
\[
 u_+(s)=\begin{pmatrix}1&s\\0&1\end{pmatrix},\qquad
 u_-(s)=\begin{pmatrix}1&0\\s&1\end{pmatrix},
\]
and set
\[
 p(x)=\frac{2-x}{x^3},\qquad q(x)=\frac{8(x-1)}{x^3},\qquad
 A(x)=u_+(p(x))u_-(q(x))=\begin{pmatrix}1+pq&p\\q&1\end{pmatrix}.
\]
Both $A$ and its inverse have no finite pole outside $\{0\}$, and $A=I+O(x^{-2})$. 
Moreover, $A(1)=u_+(1)$ and $A(2)=u_-(1)$. The Zariski closures of the cyclic subgroups generated by $A(1)$ and $A(2)$ are the two root subgroups, 
respectively, and these two root subgroups generate $\SL_2$. For example, 
the Lie algebra generated by $CE_{12}$ and $CE_{21}$ contains $E_{12}$, $E_{21}$, and $[E_{12},E_{21}]$, 
which form a basis of $\mathfrak{sl}_2$. Thus $T_A=\SL_2$, and the criterion applies without any dilation.

\subsection{A multiplicative group}\label{sub_original_8_2}
Take $B(z)=1+z^{-2}$, with singularities $0,i,-i$. Choose $t$ satisfying $t^3=2$ and consider the system
\[
 A(x)=\left(1+\frac{t^2}{x^2}\right)\in\GL_1(K),\qquad\sigma Y=A(x)Y.
\]
Since $t^{-1}\notin\Q(i)$, we have $1\notin\Span_\Q\{0,it,-it\}$, and $A=I+O(x^{-2})$. Its nonconstant coefficient has infinitely many values, which generate $\Gm$. Thus its PV group is $\Gm$. Dilation is not always necessary; for example, $1\notin\Span_\Q\{0,i,-i\}$ already holds. The chosen $t$ illustrates a uniform algebraic rule that works without first testing the individual singularity set.

\subsection{A nonsplit extension with quotient of order two}\label{sub_original_8_3}
Let
\[
 T=\{\diag(s,s^{-1}):s\in C^*\},\qquad
 h=\begin{pmatrix}0&-1\\1&0\end{pmatrix},
\]
\[ G=N_{\SL_2}(T)=T\cup hT. \]
To see this description of the normalizer, note that an element normalizing $T$ must permute its two distinct weight lines in $C^2$. It is therefore diagonal or antidiagonal; the determinant-one condition gives exactly $T\cup hT$. The quotient $G/T$ has order $2$. Conjugation by $h$ inverts $T$ and $h^2=-I$, so $(h\diag(s,s^{-1}))^2=-I$ for every $s$. Thus the extension by $C_2$ does not split: there is no element of order $2$ in the nonidentity component.

Take $t^3=2$ and $q(x)=1+t^2/x^2$. Then
\[
 A(x)=h\diag(q(x),q(x)^{-1})
 =\begin{pmatrix}0&-q(x)^{-1}\\q(x)&0\end{pmatrix}
\]
has semisimple limit $h$, no first-order term, and singularities among $0,it,-it$. Its coefficient group contains $h$ and $T$, and hence equals $G$. Its PV group is therefore the nonsplit normalizer.

\subsection{A finite cyclic group}\label{sub_original_8_4}
The construction is classical. We give a quotient-algebra presentation of the finite cyclic PV ring and verify 
its prescribed matrix representation. Let $G=\langle h\rangle$ have order $m$. 
Choose a primitive $m$-th root of unity $\zeta\in C$ and write $h=P\diag(\zeta^{e_1},\ldots,\zeta^{e_n})P^{-1}$. Put
\[
 R=K[t]/(t^m-1),\qquad\sigma(t)=\zeta t,\qquad
 Z=P\diag(t^{e_1},\ldots,t^{e_n})P^{-1}.
\]
The ring is isomorphic to $K^m$, and $\sigma$ permutes its primitive idempotents transitively. An ideal is a product of ideals in these field factors, so no nonzero proper ideal is $\sigma$-stable. Also $\sigma Z=hZ$.

Since $h$ has order $m$, $\gcd(m,e_1,\ldots,e_n)=1$. The diagonal entries of $P^{-1}ZP$ give $t^{e_i}$, and their inverses are the diagonal entries of $P^{-1}Z^{-1}P$. They belong to the algebra generated by the entries of $Z$ and $\det(Z)^{-1}$ because $Z^{-1}=\operatorname{adj}(Z)/\det(Z)$. B\'ezout's identity modulo $m$ therefore recovers $t$ in this algebra. Thus $R$ is generated by the fundamental matrix and its inverse determinant, and is a PV ring. The $K$-automorphisms are permutations of the $m$ field factors. A permutation commuting with a full cycle is determined by the image of one factor, so the commuting automorphisms are precisely $t\mapsto\zeta^j t$, with right-action matrices $h^j$.

\subsection{Omitting the singularity condition or the poles of the inverse}\label{sub_original_9_1}
Let $a(x)=1-x^{-2}$ and $y=x/(x-1)$. Then
\[
 \frac{y(x+1)}{y(x)}=\frac{(x+1)(x-1)}{x^2}=a(x).
\]
The one-dimensional system with coefficient matrix $A=(a)$ has fundamental matrix $Z=(y)\in\GL_1(K)$. Thus its PV group is trivial, while its coefficient group is $\Gm$. The normalization at infinity holds, but $\mathcal{P}(a)=\{0,1,-1\}$ violates the arithmetic condition. Counting only poles of $a$, and not of $a^{-1}$, would give the erroneous set $\{0\}$.

The degree identity (\ref{degree}) can be verified explicitly in the representation $\rho(s)=\diag(s,1)$:
\[
 v=\begin{pmatrix}x\\x-1\end{pmatrix},\qquad
 \sigma v=\frac{x}{x-1}\rho(a)v.
\]
The vector is primitive and $d=1$. The multiplier has $k_0=1$ and $k_1=-1$, so $\sum_s k_s s=-1$ and $d+\sum_s k_s s=0$.

\subsection{Omitting semisimplicity}\label{sub_original_9_2}
Let
\[
 A=\begin{pmatrix}1&1\\0&1\end{pmatrix},\qquad
 Z(x)=\begin{pmatrix}1&x\\0&1\end{pmatrix}.
\]
Then $\sigma Z=AZ$ and $Z\in\GL_2(K)$. Thus the Picard-Vessiot group is trivial, whereas $T_A\simeq\Ga$. 
The singularity set is empty, but the limit at infinity is not semisimple. 
Correspondingly, the vector $v=(x,1)^T$ satisfies $\sigma v=Av$, but the line $Kv$ has no nonzero constant generator. This is precisely the failure of $\ker(M-I)\cap\im(M-I)=0$ in the proof of rigidity.

\subsection{Allowing a first-order term at infinity}\label{sub_original_9_3}
Consider
\[
 A(x)=\begin{pmatrix}1+x^{-1}+x^{-2}&-x^{-1}\\x^{-2}&1-x^{-1}\end{pmatrix},
 \qquad v(x)=\begin{pmatrix}x\\1\end{pmatrix}.
\]
A direct computation shows that $\det A=1$ and $\sigma v=Av$. 
Both $A$ and $A^{-1}$ have finite poles only at $0$, and $A(\infty)=I$ is semisimple. Nevertheless, 
$Kv$ is not a constant line. Thus the vanishing of the $x^{-1}$ coefficient cannot be omitted from the hypotheses.

More generally, if $B=H+B_1/x+O(x^{-2})$ and $\sigma v=cBv$, the same two-coefficient calculation gives, for $M=\gamma H$,
\[
 Mv_d=v_d,\qquad(d+q)v_d=(M-I)v_{d-1}+\gamma B_1v_d.
\]
When $H$ is semisimple, projection $\Pi_1$ onto the $1$-eigenspace of $M$ gives $(d+q)v_d=\gamma\Pi_1B_1v_d$. The proof requires this extra term to vanish; our construction arranges the stronger condition $B_1=0$ uniformly in every representation.

\subsection{Rational gauges of constant systems}\label{sub_original_9_5}
\begin{lemma}[A scalar coboundary creates an integer resonance]\label{result_original_9_1}
If $f\in C(x)^*$ is nonconstant, the finite zero--pole support of $\sigma(f)/f$ contains two distinct points differing by a nonzero integer.
\end{lemma}
\begin{proof}
A nonconstant rational function has a finite zero or pole. In an orbit $a+\Z$ containing one, put $n_k=\ord_{a+k}(f)$. This is a nonzero finitely supported sequence. Let $m$ and $M$ be the least and greatest indices with $n_k\ne0$. Then
\[ \ord_{a+k}(\sigma(f)/f)=n_{k+1}-n_k \]
is nonzero at $k=m-1$ and at $k=M$. These are distinct indices, and the corresponding points differ by $M-m+1\in\Z_{>0}$.
\end{proof}

Suppose $A$ satisfies \eqref{hyp} and $A=\sigma(P)DP^{-1}$ with $P\in\GL_n(K)$ and constant semisimple $D$. Replacing $(P,D)$ by $(PS,S^{-1}DS)$ for a suitable $S\in\GL_n(C)$, we may assume that $D=\diag(\lambda_1,\ldots,\lambda_n)$. Applying Theorem~\ref{rigidity} to the columns of $P$ gives
\[
 P=Q\diag(f_1,\ldots,f_n),\qquad Q\in\GL_n(C),\quad f_j\in K^*.
\]
Consequently,
\[
 Q^{-1}AQ=\diag\left(\lambda_1\frac{\sigma(f_1)}{f_1},\ldots,
 \lambda_n\frac{\sigma(f_n)}{f_n}\right).
\]
The finite zeros and poles of each $\sigma(f_j)/f_j$ lie in $\mathcal P(A)$. Lemma~\ref{result_original_9_1} and $1\notin\Span_\Q\mathcal P(A)$ therefore force every $f_j$ to be constant. Hence $P$ is constant.

\section{An alternative interpolation construction}\label{sec_original_10}
In this section, we give an alternative construction of the generating rational curve using a rational chart 
and Hermite interpolation. The interpolation prescribes tangent vectors spanning the Lie algebra of the identity 
component and ensures the required expansion at infinity. Smoothness and connectedness then imply that the values of the curve generate the identity component.

A connected linear algebraic group $H$ over $C$ is a rational variety; see \cite[proof of Proposition 3.3(2)]{PS}. 
Put $d=\dim H$. Translating a rational chart gives an open subset $U\subset\mathbb A^d$ and a regular map that is an isomorphism onto an open subset,
\[
 F:U\longrightarrow H,\qquad0\in U,\quad F(0)=I,\quad
 dF_0:C^d\xrightarrow{\sim}\Lie(H).
\]
For $d>0$, put $b_i=i$, $M=2d+1$, and
\[
 L_i(z)=\prod_{j\ne i}\frac{z-b_j}{b_i-b_j},\qquad
 u_i(z)=\frac{b_i^M(z-b_i)L_i(z)^2}{z^M}.
\]
At $z=b_i$, the factor $z-b_i$ gives derivative $1$ after multiplication by $b_i^M/z^M$. At $z=b_j$, $j\ne i$, the square of $L_i$ has a double zero. Hence
\[
 u_i(b_j)=0,\qquad u_i'(b_j)=\delta_{ij},\qquad u_i=O(z^{-2}).
\]
The last estimate follows because the numerator has degree $2d-1$. The curve $B(z)=F(u_1(z),\ldots,u_d(z))$ is 
defined generically and at the specified nodes and infinity. 
At each interpolation node and at infinity, the parameter vector $(u_1(z),\ldots,u_d(z))$ equals $0\in U$. Hence the composition defining $B$ is well defined at these points. It satisfies
\[
 B(z)=I+O(z^{-2}),\qquad B(b_j)=I,\qquad B'(b_j)=dF_0(e_j).
\]
If a closed subgroup $P\subseteq H$ contains all regular values of $B$, its Lie algebra contains these derivatives, which form a basis of $\Lie(H)$. Smoothness in characteristic zero gives $\dim P=\dim H$; connectedness of $H$ implies $P=H$. Its finite singularities need not be algebraic over $\Q$, but the first assertion of Lemma~\ref{result_original_7_1} works for every finite set. Thus $A(x)=hB(x/t)$ gives a second construction under the same rigidity criterion.
\section{Acknowledgements}
The author was supported by the Strategic Priority Research Program of the Chinese Academy of Sciences (Grant Nos. XDA0480502 and XDA0480503) and by the National Natural Science Foundation of China (Grant No. 12471479).
The author is very grateful to Michael Singer for his valuable suggestions concerning the results of this paper.

\end{document}